\documentclass{article}
\usepackage{amsmath,amsfonts}

\DeclareMathSymbol\nullset{\mathord}{AMSb}{"3F}

\begin{document}
\newtheorem{theorem}{Theorem}[section]
\newtheorem{corollary}{Corollary}[section]
\newtheorem{lemma}{Lemma}[section]
\newtheorem{proposition}{Proposition}[section]
{\newtheorem{definition}{Definition}[section]}

\title{The identities and the central polynomials of the infinite dimensional 
             unitary Grassmann algebra over a finite field}
	     
\author{C{.} Bekh-Ochir and S{.} A{.} Rankin
  }

\maketitle

\begin{abstract}
We describe the $T$-ideal of identities and the $T$-space of central polynomials for the 
infinite dimensional unitary Grassmann algebra over a finite field.
\end{abstract}



\newcounter{parts}
\def\set#1\endset{\{\,#1\,\}}
\def\gen#1{\left<{#1}\right>}
\def\rest#1{{}_{{}_{#1}}}
\def\com#1,#2{[{#1},{#2}]}
\def\choice#1,#2{\binom{#1}{#2}}
\def\kx{k\langle X\rangle}
\def\kzerox{k_0\langle X\rangle}
\def\konex{k_1\langle X\rangle}
\def\unitgrass{G}
\def\nonunitgrass{G_0}
\def\finiteunitgrass#1{G(#1)}
\def\finitenonunitgrass#1{G_0(#1)}
\def\siderovset{SS}
\def\boundedsiderovset{BSS}
\def\finitesiderovsetone#1{SS'(#1)}
\def\finitesiderovsettwo#1{SS''(#1)}
\def\boundedsiderovset{BSS}
\def\bss#1{BSS(#1)}
\def\lend(#1){lend(#1)}
\def\lbeg#1{lbeg(#1)}
\def\endof(#1){end(#1)}
\def\beg(#1){beg(#1)}
\def\finitesidset#1{SS(#1)}
\let\cong=\equiv
\def\mod#1{\,\,(\text{mod}\,#1)}
\def\mz(#1){M_0(#1)}
\def\wt(#1){\text{wt}(#1)}
\def\dom(#1){\text{dom}(#1)}
\def\from{\mkern2mu{:}\mkern2mu}

\def\proof{\ifdim\lastskip<\smallskipamount\relax\removelastskip
  \vskip\smallskipamount\fi\leavevmode\hbox to 10pt{\hfil}{\it Proof. }}
\def\strutdepth{\dp\strutbox}
\def\epmarker{\vbox to \strutdepth{\baselineskip\strutdepth\vss\hfill{%
\hbox to 0pt{\hss\vrule height 4pt width 4pt depth 0pt}\null}}}
\def\edproofmarker{\strut\vadjust{\kern-2\strutdepth\epmarker}}
\def\endproof{\edproofmarker\vskip10pt}

\section{Introduction}
 In 1987, one of the fundamental results in the theory of PI-algebras
 was obtained by A\hbox{.} R\hbox{.} Kemer (\cite{Ke}). Kemer proved
 that every system of identities in an associative algebra over a field
 of characteristic zero is finitely based, which provided a positive
 answer to a question raised by W\hbox{.} Specht (\cite{Sp}) in 1950.
 Shortly after Kemer's result appeared,  A\hbox{.} V\hbox{.} Grishin
 introduced the concept of $T$-space (\cite{Gr1}, \cite{Gr2}); a vector
 subspace of an algebra that is closed under the natural action of the
 monoid $T$ of all endomorphisms of the algebra. As shown by Grishin and
 V\hbox{.} V\hbox{.} Shchigolev in the influential  survey paper
 \cite{GrSh}, $T$-spaces have important applications in the theory of
 PI-algebras and in the problem of the finite-basedness of $T$-ideals.

 $T$-ideals arise in the study of the identities of an associative
 algebra, and very closely related to the $T$-ideal of identities of an
 associative algebra is the $T$-space of central polynomials of an
 associative algebra; the set of all elements that map into the centre 
 under every algebra homomorphism from the free associative algebra into the
 given associative algebra. A\hbox{.} Y\hbox{.} Belov, writing in
 \cite{Be}, observed that if one regards PI-theory as a kind of
 viewpoint for non-commutative algebraic geometry, then the Grassmann
 algebra serves as one of the most important examples of new objects
 that are analogues of prime algebras. In this context, it seemed
 natural to investigate the $T$-space of central polynomials of the
 Grassmann algebra, with a view to determining whether or not this
 $T$-space is Spechtian.
  
 In \cite{Ra} and \cite {CR}, we identified the $T$-space of central 
 polynomials of the finite and the infinite dimensional, unitary and
 nonunitary Grassmann algebras over a field of arbitrary characteristic, 
 although only for an infinite field in the case of the unitary Grassmann
 algebra.  In these earlier works, we were able to utilize descriptions
 of the $T$-ideal of identities for the corresponding Grassmann algebras
 due to Chiripov and Siderov \cite{Si}, Giambruno and Koshlukov
 \cite{Gi}, and Stojanova-Venkova \cite{At}, but for the unitary
 Grassmann algebras over a finite field, the $T$-ideal of identities was
 not yet known. A\hbox{.} Regev had initiated a study of the identities
 of the infinite dimensional unitary Grassmann algebra over a finite
 field in \cite{Re}, but a complete description of the identities for
 that case was not forthcoming.
 
 The purpose of this paper is to present a complete description  of the
 $T$-ideal of the identities of the infinite dimensional unitary
 Grassmann algebra over a finite field,  thus completing work that was
 begun by Regev \cite{Re} in 1991. We then use this information to
 provide a complete description of the $T$-space of central polynomials
 in this case as well.

 In a subsequent article (see \cite{CR5}), we establish that if $p>2$
 and $k$ is an arbitrary field of  characteristic $p$, then neither the
 $T$-space of central polynomials of the unitary nor the nonunitary
 infinite dimensionalGrassmann algebra over $k$ is finitely based.
 
\section{Preliminaries}
 Let $k$ be a finite field of characteristic $p$ and size $q$, and let
 $X$ be a countably infinite set,  say $X=\set x_i\mid i\ge 1\endset$.
 Then $\kzerox$ denotes the free (nonunitary) associative  $k$-algebra
 over $X$, while $\konex$ denotes the free unitary associative
 $k$-algebra over $X$. 

 Let $A$ denote any associative $k$-algebra.  Any linear subspace of $A$
 that is invariant under the natural action of the monoid $T$ of all
 algebra endomorphisms of $A$ is called a $T$-space of $A$, and if a
 $T$-space happens  to also be an ideal of $H$, then it is called a
 $T$-ideal of $A$. For $B\subseteq A$, the  smallest $T$-space
 containing $B$ shall be denoted by $B^S$, while the smallest $T$-ideal 
 of $A$ that contains $B$ shall be denoted by $B^T$. In this article, we
 shall deal only with $T$-spaces and $T$-ideals of $\kzerox$ and
 $\konex$.

 A nonzero element $f\in \kzerox$ is called an {\em identity} of $A$ if
 $f$ is in the kernel of every $k$-algebra homomorphism from $\kzerox$
 to $A$ (every unitary $k$-algebra homomorphism from $\konex$ if $A$ is
 unitary).  The set consisting of $0$ and all identities of $A$ is a
 $T$-ideal of $\kzerox$ (and of $\konex$ if $A$ is unitary),  denoted by
 $T(A)$. An element $f\in \kzerox$ is called a {\it central polynomial}
 of $A$  if $f\notin T(A)$ and the image of $f$ under any $k$-algebra
 homomorphism from $\kzerox$ (unitary $k$-algebra homomorphism from
 $\konex$ if $H$ is unitary) to $A$ belongs to $C_{A}$, the centre of
 $A$. The $T$-space of $\kzerox$ (or of $\konex$ if $A$ is unitary) that
 is generated by the set of all central polynomials of $A$ is denoted by
 $CP(A)$.

 Let $\unitgrass$ denote the (countably) infinite dimensional unitary
 Grassmann algebra over $k$, so there exist $e_i\in \unitgrass$, $i\ge
 1$,  such that  for all $i$ and $j$, $e_ie_j=-e_je_i$, $e_i^2=0$, and
 $\mathcal{B}=\set e_{i_1}e_{i_2}\cdots e_{i_n}\mid n\ge 1,
 i_1<i_2<\cdots i_n\endset$, together with $1$, forms a linear basis for
 $G$. Let $E$ denote the set $\set e_i\mid i\ge1\endset$. The subalgebra
 of $\unitgrass$ with linear basis $\mathcal{B}$ is the infinite
 dimensional nonunitary Grassmann algebra over $k$, and is denoted by
 $\nonunitgrass$. Then for any positive integer $m$, the unitary
 subalgebra of $\unitgrass$ that is generated by $\set
 e_1,e_2,\ldots,e_m\endset$, is denoted by $\finiteunitgrass{m}$, while
 the nonunitary subalgebra of $\nonunitgrass$ that is generated by the
 same set is denoted by $\finitenonunitgrass{m}$.

 It is well known that $T^{(3)}$, the $T$-ideal of $\konex$ 
 generated  by $\com {\com x_1,{x_2}},{x_3}$, is contained in
 $T(\unitgrass)$. For convenience,  we shall write $\com x_1,{x_2,x_3}$
 for  $\com {\com x_1,{x_2}},{x_3}$.

 In that paper, Regev showed that $\set
 [x,y,z],x^{qp}-x^p\endset^T\subseteq T(\unitgrass)$. By working modulo
 the $T$-ideal $\set [x,y,z],x^{qp}-x^p\endset^T$, we are able to
 establish that in fact, equality holds when $p>2$, while for $p=2$,
 $T(\unitgrass)=\set \com x_1,{x_2},x_1^{2q}-x^2\endset^T$.  Then, with
 full knowledge of $T(\unitgrass)$ in hand, we are able to obtain the
 $T$-space of central polynomials of the infinite dimensional unitary
 Grassmann algebra over a finite field (the one outstanding case).

 Evidently (since all Grassmann algebras over a field of characteristic
 2 are commutative), $CP(\unitgrass)=\konex$ if $p=2$, and we show that
 for $p>2$, 
 \begin{align*}
   CP(\unitgrass)&=T(\unitgrass)+\set \com x_1,{x_2},x_1^p\endset^S+\set x_{1}^p
    \prod_{i=1}^k\com x_{2i},{x_{2i+1}}x_{2i}^{p-1}x_{2i+1}^{p-1}\mid k\ge1\endset^S\\
    &=\set \com x_1,{x_2,x_3}{x_4},x_1(x_2^{qp}-x_2^p)\endset^S+\set\com x_1,{x_2},x_1^p\endset^S\\
    &\hskip80pt+\set x_{1}^p
    \prod_{i=1}^k\com x_{2i},{x_{2i+1}}x_{2i}^{p-1}x_{2i+1}^{p-1}\mid k\ge1\endset\bigr)^S.
\end{align*}    

 We complete this section with a brief description of results from the
 literature that will be required in this work. To begin with, the
 following lemma summarizes  discussion found in Chiripov and Siderov
 \cite{Si}. A product term $e_{i_1}e_{i_2}\cdots e_{i_n}$ in
 $\nonunitgrass$ is said to be {\it even} if $n$ is even, otherwise the
 product term is said to be {\it odd}. $u\in \nonunitgrass$ is said to
 be {\it even} if $u$ is a linear combination  of even product terms,
 while $u$ is said to be {\it odd} if $u$ is a linear combination of odd
 product terms. Let $C$ denote the set of all even elements of
 $\nonunitgrass$, and let $H$ denote the set of all odd elements of
 $\nonunitgrass$. Note that $C$ and $H$ are $k$-linear subspaces of
 $\nonunitgrass$, and $C$ is closed under multiplication, $H^2\subseteq
 C$, and $CH=HC\subseteq H$. Evidently, $\nonunitgrass=C\oplus H$ as
 $k$-vector spaces.
 
\begin{lemma}\label{lemma: useful}\ 
   \vskip-1.75\baselineskip\null  
   \begin{list}{(\roman{parts})}{\usecounter{parts}}
  \item 
    $C_{\nonunitgrass}=C$, and $C_{\unitgrass}=k\oplus C$.
  \item 
    For $h,u\in H$, $hu=-uh$. In particular, $h^2=0$ (since $p>2$).
  \item 
    Let $g\in \nonunitgrass$, so there exist (unique) $c\in C$ and $h\in H$ such that $g=c+h$. For any 
  positive integer $n$, $g^n=c^n+nc^{n-1}h$. 
  \item 
   For $g\in \nonunitgrass$, $g^p=0$, and for any $\alpha\in k$, $(\alpha+g)^p=\alpha^p$.
  \item 
   Let $c_1,c_2\in C$ and $h_1,h_2\in H$, and set $g_1=c_1+h_1$, $g_2=c_2+h_2$. Then
   for any nonnegative integers $m_1,m_2$, $\com g_1,{g_2}g_1^{m_1}g_2^{m_2}=2c_1^{m_1}c_2^{m_2}h_1h_2$
   (where $g_i^0$ and $c_i^0$ are understood to mean that the factors $g_i^0$ and $c_i^0$ are omitted).
  \item 
   Let $u\in \nonunitgrass$. Then $u^{n+1}=0$, where $n$ is the number of distinct basic product terms in the expression
   for $u$ as a linear combination of elements of $\mathcal{B}$.
 \end{list}
\end{lemma}

\begin{definition}\label{definition: siderov's elements}
 Let $\siderovset$ denote the set of all elements of the form 
 \begin{list}{(\roman{parts})}{\usecounter{parts}}
  \item $\prod_{r=1}^t x_{i_r}^{\alpha_r}$, or
  \item $\prod_{r=1}^s \com x_{j_{2r-1}},{x_{2r}}x_{j_{2r-1}}^{\beta_{2r-1}}x_{j_{2r}}^{\beta_{2r}}$, or
  \item $\bigl(\prod_{r=1}^t x_{i_r}^{\alpha_r}\bigr)\prod_{r=1}^s \com x_{j_{2r-1}},{x_{2r}}x_{j_{2r-1}}^{\beta_{2r-1}}x_{j_{2r}}^{\beta_{2r}}$, 
 \end{list}  
 \noindent where $j_1<j_2<\cdots j_{2s}$, $\beta_i\ge0$ for all $i$, $i_1<i_2<\cdots < i_t$, 
 $\set i_1,\ldots,i_r\endset\cap\set j_1,\ldots,j_{2s}\endset=\nullset$, and $\alpha_i\ge 1$ for all $i$.

 Let $u\in \siderovset$. If $u$ is of the form (i), then the beginning
 of $u$, $\beg(u)$, is $\prod_{r=1}^t x_{i_r}^{\alpha_r}$, the end of
 $u$, $\endof(u)$, is empty, the length of the beginning of $u$,
 $\lbeg{u}$, is equal to $t$ and the length of the end of $u$,
 $\lend(u)$, is 0. If $u$ is of the form (ii), then we say that
 $\beg(u)$, the beginning of $u$, is empty, $\endof(u)$, the end of $u$,
 is $\prod_{r=1}^s \com
 x_{j_{2r-1}},{x_{2r}}x_{j_{2r-1}}^{\beta_{2r-1}}x_{j_{2r}}^{\beta_{2r}}$,
 and $\lbeg{u}=0$ and $\lend(u)=s$. If $u$ is of the form (iii),  then
 we say that $\beg(u)$, the beginning of $u$, is $\prod_{r=1}^t
 x_{i_r}^{\alpha_r}$, $\endof(u)$, the end of $u$, is $\prod_{r=1}^s
 \com
 x_{j_{2r-1}},{x_{2r}}x_{j_{2r-1}}^{\beta_{2r-1}}x_{j_{2r}}^{\beta_{2r}}$,
 and $\lbeg{u}=t$ and $\lend(u)=s$.
\end{definition}

In \cite{Si},  Siderov introduced a total order on the set $\siderovset$
which was useful in his work on the identities of the infinite
dimensional nonunitary Grassmann algebra. 

\begin{definition}[Siderov's ordering]\label{definition: total order}
 For $u,v\in \siderovset$, we say that $u>v$ if one of the following requirements holds.
  \begin{list}{(\roman{parts})}{\usecounter{parts}}
  \item $\deg u > \deg v$.
  \item $\deg u = \deg v$ but $\lend(u)<\lend(v)$.
  \item $\deg u = \deg v$ and $\lend(u)=\lend(v)$, but there exists $i \ge 1$ such that 
        $\deg_{x_i} u > \deg_{x_i} v$ and for each $j<i$, $\deg_{x_j} u = \deg_{x_j} v$.
  \item $\deg u = \deg v$, $\lend(u)=\lend(v)$ and for each $i \ge 1$,
        $\deg_{x_i} u = \deg_{x_i} v$, and there exists $j \ge 1$ such that 
        $x_j$ appears in $\beg(u)$ and in $\endof(v)$, and 
        for each $k>j$, $x_k$ appears in $\beg(u)$ if and only if $x_k$ appears in 
        $\beg(v)$.
  \end{list}
 \end{definition}


\section{The $T$-ideal of identities of the infinite dimensional unitary Grassmann algebra over a finite field}

\begin{lemma}[\cite{Re}, Lemma 1.5]\label{lemma: regev's identities} 
 $x^{qp}-x^p\in T(\unitgrass)$. Moreover, if $f(x)$ is a one-variable identity
 of $\unitgrass$, then $x^{qp}-x^p$ divides $f$ in $\konex$.
\end{lemma}

\begin{definition}\label{bounded siderov elements}
  Let
  \begin{align}
   \boundedsiderovset&=\{\, u\in \siderovset\mid \text{for each $i$, }\deg_{x_i}(u)<p\text{ if $x_i$ appears in $\beg(u)$} \notag\\
      &\hskip110pt\text{or $\deg_{x_i}(u)\le p$ if $x_i$ appears in $\endof(u)$}\,\}.\notag
  \end{align}
\end{definition}
  
\begin{definition}\label{definition: support}  For
$u=e_{i_1}e_{i_2}\cdots e_{i_n}\in \mathcal{B}$, let $s(u)=\set
e_{i_1},e_{i_2},\ldots,e_{i_n}\endset$ and  $\wt(u)=|s(u)|$, while
$s(1)=\nullset$ and $\wt(1)=0$. We call $s(u)$ and $\wt(u)$ the support
and weight of $u$,   respectively. Now for any $g\in \unitgrass$, $g\ne
0$, $g=\sum_{i=1}^m a_ig_i$ with $a_i\in k^*=k-\set 0\endset$ and
$g_i\in \mathcal{B}\cup\set1\endset$.   Let $s(g)=\bigcup_{i=1}^m
s(g_i)$, $\wt(g)=\max\set\wt(g_i)\mid i=1,2,\ldots,m\endset$, and  
$\dom(g)=\sum_{\wt(g_i)=\wt(g)} a_ig_i$, while we define $s(0)=\nullset$
and $\wt(0)=0$.   We call $s(g)$ the support of $g$, $\wt(g)$ the weight
of $g$, and  $\dom(g)$ the dominant part of $g$. Note that if
$s(g_1)\cap s(g_2)=\nullset$, then $\dom(g_1g_2)=\dom(g_1)\dom(g_2)$ 
and $\wt(g_1g_2)=\wt(g_1)+\wt(g_2)$.
\end{definition}

\begin{lemma}\label{technical result}
 Let $n$ and $\gamma$  be positive integers and let $\lambda\in k$. Then the following hold.
 \begin{list}{(\roman{parts})}{\usecounter{parts}}
 \item
  $\dom({(\lambda+\sum_{\epsilon=1}^n e_{2\epsilon-1}e_{2\epsilon})^\gamma})
 = \begin{cases}
 \frac{\gamma!}{(\gamma-n)!}\lambda^{\gamma-n}\prod_{\epsilon=1}^{2n}e_{\epsilon} &\text{if $\gamma\ge n$}\\
 \gamma!\sum_{\substack{ J\subseteq J_n\\|J|=\gamma}}\prod_{j\in J} e_{2j-1}e_{2j} &\text{if $\gamma<n$.}
 \end{cases}
 $
 \item
 \begin{align*}
  \dom({(\lambda+e_{2n+1}+\sum_{\epsilon=1}^n e_{2\epsilon-1}&e_{2\epsilon})^\gamma})
 =\\
 &\begin{cases}
 \frac{\gamma!}{(\gamma-n)!}\lambda^{\gamma-n}\prod_{\epsilon=1}^{2n+1}e_{\epsilon} &\text{if $\gamma> n$}\\
 \gamma!\sum_{\substack{ J\subseteq J_n\\|J|=\gamma}}\prod_{j\in J} e_{2j-1}e_{2j} &\text{if $\gamma\le n$.}
 \end{cases}
 \end{align*}
 \end{list}
\end{lemma}

\begin{proof}
 For (i), note that $\dom(\lambda+\sum_{\epsilon=1}^n e_{2\epsilon-1}e_{2\epsilon})^\gamma$ is equal to
 $$
  \dom(\sum_{\substack{c_0+c_1+\cdots+c_{n}=\gamma\\
  c_0,\,c_1,\ldots,\,c_n\ge0}} \choice \gamma,{c_0,\,c_1,\,\ldots,\,c_n}\lambda^{c_0}
 \prod_{\epsilon=1}^n (e_{2\epsilon-1}e_{2\epsilon})^{c_\epsilon}).
 $$
 The dominant part will therefore be obtained by setting as many as
 possible of the $c_i$ values to $1$, $i\ge1$ (if $c_i\ge2$ when
 $i\ge1$, the summand will be 0). The  result is as shown in (i).

 For (ii), note that $\dom({\lambda+e_{2n+1}+\sum_{\epsilon=1}^n e_{2\epsilon-1}e_{2\epsilon})^\gamma})$ is equal to
 $$
 \dom({\sum_{\substack{c_0+c_{n+1}+c_1+\cdots+c_{n}=\gamma\\ c_0,\,c_{n+1},\,c_1,\ldots,\,c_n\ge0}}
   \choice \gamma,{c_0,\,c_{n+1},\,\,c_1,\,\ldots,\,c_n}\lambda^{c_0}
   e_{2n+1}^{c_{n+1}}\prod_{\epsilon=1}^n (e_{2\epsilon-1}e_{2\epsilon})^{c_\epsilon}}).
 $$
 The dominant part will therefore be obtained
 by setting as many as possible of the $c_i$ values to $1$, $i\ge1$ (if $c_i\ge2$ when $i\ge1$, the summand will be 0). The
 result is as shown in (ii).
\end{proof}

 Let $proj_k\from \unitgrass\to k$ denote the $k$-algebra homomorphism
 that is determined by mapping $1$ to $1$, and $e_i$ to $0$.
 
\begin{proposition}\label{proposition: main evaluation fact}
 Let $u\in\boundedsiderovset$, and set $m=2\deg(u)-2\lend(u)$. For each $i\ge1$, let $\lambda_i\in k$.
 Then there exists a homomorphism $\varphi\from \konex\to
 \finiteunitgrass{m}$ such that the following hold.
   \begin{list}{(\roman{parts})}{\usecounter{parts}}
   \item For each index $i$, $proj_k(\varphi(x_i))=\lambda_i$.
  \item $\displaystyle\dom(\varphi(u))= 2^{\lend(u)}\mkern-25mu\prod_{x\text{ in }\beg(u)}\mkern-25mu \deg_{x}(u)\mkern2mu!
    \mkern-25mu\prod_{x\text{ in }\endof(u)}\mkern-25mu (\deg_{x}(u)-1)\mkern2mu!\mkern5mu\prod_{i=1}^{m}e_i$.
  \item For any $v\in\boundedsiderovset$ with $u>v$, $m=\wt(\varphi(u))>\wt(\varphi(v))$.
  \end{list}
\end{proposition}

\begin{proof}
 The homomorphism $\varphi$ is determined by the following assignments.
 First, any variable $x_i\in X$ that does not appear in $u$ is mapped to
 $\lambda_i$. Then for any variable $x$ that appears in $\beg(u)$,
 choose an index offset $N=N_x$, $E_x=\set e_{N+\epsilon}\mid
 \epsilon=1,2,\ldots,2\deg_{x}(u)\endset\subseteq E$, and map $x$ to
 $\lambda_x+\sum_{\epsilon=1}^{\deg_{x}(u)}
 e_{N+2\epsilon-1}e_{N+2\epsilon}$. Finally, for any variable $x$ that
 appears in $\endof(u)$, choose an index offset $N=N_x$, $E_x=\set
 e_{N+\epsilon}\mid \epsilon=1,2,\ldots,2\deg_{x}(u)-1\endset\subseteq
 E$, and map $x$ to
 $\lambda_x+e_{N+2\deg_{x}(u)-1}+\sum_{\epsilon=1}^{\deg_{x}(u)-1}
 e_{N+2\epsilon-1}e_{N+2\epsilon}$. Note that (i) is satisfied by this
 assignment. The offsets $N_x$ are chosen so that $x\ne y$ implies that 
 $E_x\cap E_y=\nullset$ and $\bigcup_{x\text{ appears in }u} E_x=\set
 e_i\mid i=1,2,\ldots,m\endset$.
 
 Recall that for $g_1,g_2\in\unitgrass$, $\dom(g_1g_2)=\dom(g_1)\dom(g_2)$
 if $s(g_1)\cap s(g_2)=\nullset$. In
 particular, since $u\in\boundedsiderovset$ (where the cases of $u$ with
 empty beginning or empty end are just simplifications of the following
 argument), $u$ is of the form  
 $$
  \prod_{r=1}^tx_{i_r}^{\alpha_r}\prod_{r=1}^s [\mkern1mu x_{j_{2r-1}},
  {x_{j_{2r}}}]x_{j_{2r-1}}^{\beta_{2r-1}}x_{j_{2r}}^{\beta_{2r}},
 $$
 where for each $r=1,2,\ldots,t$, $1\le \alpha_{r}\le p-1$ and for each $r=1,2,\ldots,2s$, $0\le \beta_{r}\le p$, so
 $\dom(\varphi(u))$ will be the product of $\dom(\varphi(x_{i_r})^{\alpha_r}))$, $r=1,2,\ldots,t$, and 
 $$
  \dom([\mkern1mu \varphi(x_{j_{2r-1}}),\varphi({x_{j_{2r}}})]\varphi(x_{j_{2r-1}})^{\beta_{2r-1}}\varphi(x_{j_{2r}})^{\beta_{2r}})
 $$
 for $r=1,2,\ldots,s$. 

 We now apply Lemma \ref{technical result} to evaluate the dominant part
 of $\varphi(u)$, where for convenience, we shall let $N_{i}$ denote
 $N_x$ where $x=x_{i}$. Note that for $g=\lambda+c+h$, where $\lambda\in
 k$, $c\in C$ and $h\in H$, we have for any $g_1\in\unitgrass$ that
 $\com g,{g_1}g^r=\com h,{g_1}(\lambda+c)^r$. If $s(g)\cap
 s(g_1)=\nullset$, then $g$'s contribution to the dominant part is
 $\dom(h(\lambda+c)^r)$. We now apply Lemma \ref{technical result} to
 obtain 
 \begin{align*}
 \dom(\varphi(u))&=\prod_{r=1}^t 
  \biggl(\alpha_r!\prod_{\epsilon=1}^{2\alpha_{r}}e_{N_{i_r}+\epsilon}\biggr)\mkern4mu 2^{\lend(u)}\\
  &\hskip50pt\prod_{r=1}^s  \biggl( \beta_{2r-1}!
  \prod_{\epsilon=1}^{2\beta_{2r-1}+1} \mkern-20mu e_{N_{j_{2r-1}}+\epsilon}\mkern6mu  
  \beta_{2r}!\,\prod_{\epsilon=1}^{2\beta_{2r}+1}e_{N_{j_{2r}}+\epsilon}\biggr) \\
  &=2^{\lend(u)}\prod_{r=1}^t\alpha_r!\prod_{r=1}^{2s}\beta_{r}!\mkern20mu \prod_{\epsilon=1}^m e_\epsilon
 \end{align*} 
 and so (ii) holds. 
 
 Finally, suppose that $v\in \boundedsiderovset$, and that $u>v$. Note
 that $\wt(\varphi(u))=m>0$. If a variable appears in $\beg(u)$ and in
 $\endof(v)$, then $\varphi(v)=0$ and so $\wt(\varphi(v))=0$. Thus we
 may assume that every variable that appears in $\beg(u)$ and in $v$
 then appears in $\beg(v)$. Since conditions (ii) and (iv) of the
 definition of the Siderov ordering imply that there exists a variable
 that appears in $\beg(u)$ and in $\endof(v)$, we see that $u>v$ must
 hold by virtue of conditions (i) or (iii) of the definition. We now
 calculate the weight of the dominant part of $\varphi(v)$. First,
 observe that the contribution to the weight of the dominant part of
 $\varphi(v)$ of $x\in X$ that appears in $v$ but  not in $u$ is 0.
 Next, by Lemma \ref{technical result}, the weight of the contribution
 of $x\in X$ that appears in both $u$ and $v$  is $\min\set
 2\deg_{x}(u),2\deg_x(v)\endset$ if $x$ appears in $\beg(u)$ (so the
 contribution of $x$ to the dominant part of $\varphi(u)$ has weight
 $2\deg_x(u)$), while it is $\min\set 2\deg_x(u)-1, 2\deg_x(v)-1\endset$
 if  $x$ appears in $\endof(u)$ and $\endof(v)$ (so the contribution of
 $x$ to the dominant part of $\varphi(u)$ has weight $2\deg_x(u)-1$),
 and it is $\min\set 2\deg_x(u)-1, 2\deg_x(v)\endset$ if $x$ appears in
 $\endof(u)$ and in $\beg(v)$ (so the contribution of $x$ to the
 dominant part of $\varphi(u)$ has weight $2\deg_x(u)$). Now, either of
 conditions (i) or (iii) implies that there is a variable $x$ such that
 $\deg_{x}(u)>\deg_{x}(v)\ge0$, so it follows that
 $\wt(\varphi(u))>\wt(\varphi(v))$. 
\end{proof}

The following result is well known.

\begin{lemma}\label{lemma: p powers commute}
 $\com x_1^p,{x_2}\in T^{(3)}$.
\end{lemma}

\begin{lemma}\label{lemma: field identities}
Let $f\in \konex$ be of the form $f=\sum_{r=1}^t \lambda_r u_r$, where for each $r=1,2,\ldots,t$,
$u_r\in \siderovset$, $\lend(u_r)=0$, and for each $x\in X$ that appears in $u_r$, $\deg_x(u_r)<q$. 
If $f\in T(k)$, then $f=0$.
\end{lemma}

\begin{proof}
 The proof will be by induction on $n$, the number of variables that
 appear in $f$.  If $f\in T(k)$ is a single variable polynomial, then
 $f$ is divisible by $x^q-x$, and thus $f=0$ by degree considerations.
 Suppose that $n\ge1$ is such that the assertion holds and consider an
 $n+1$  variable polynomial $f\in T(k)$ of the required form. Suppose
 that $f\ne 0$. We may assume that the variables of $f$ are
 $x_1,x_2,\ldots,x_{n+1}$.  Let $m=\deg_{x_{n+1}}(f)$.  Then we may
 write $f=\sum_{r=1}^m f_rx_{n+1}^{r}$, where for each $r$,  $f_r$ is a
 linear combination of elements of $\siderovset$ with empty end in which
 each variable has  degree less than $q$, but on only the $n$ variables
 $x_1,x_2,\ldots,x_{n}$, and $f_m\ne 0$. It follows from the induction
 hypothesis that $f_m$ is not an  identity, and so there are
 $g_1,g_2,\ldots,g_n\in k$ such that $f_m(g_1,\ldots,g_n)\ne 0$. But
 then $f(g_1,g_2,\ldots,g_n,x)$  is a one variable identity of degree
 $m<q$, which implies that $f(g_1,g_2,\ldots,g_n,x)$ is the zero
 polynomial. In particular, $f_m(g_1,\ldots,g_n)= 0$, which is a
 contradiction. Thus $f=0$, as required, and the result follows by
 induction.
\end{proof} 

\begin{definition}\label{def of p-polynomial}
 $f\in \konex$ shall be called a {\em $p$-polynomial} if either $f\in k$
 or else $f=\sum_{r=1}^t \lambda_r u_r$, where for each 
 $r=1,2,\ldots,t$, $u_r\in \siderovset$, $\lend(u_r)=0$, and for each
 $x\in X$ that appears in $u_r$,  $\deg_x(u_r)<qp$ and $\deg_x(u_r)\cong
 0\mod{p}$.
\end{definition}

\begin{corollary}\label{corollary: needed for main identity theorem}
 If $f\in T(\unitgrass)$ is a $p$-polynomial, then $f=0$.
\end{corollary}

\begin{proof}
 If $f\in k$, then $f$ is obviously $0$, so consider $f\notin k$.
 Let $v_i\in \siderovset$ be such that $u_i(x_1,x_2,\ldots,x_m)=v_i(x_1^p,x_2^p,\ldots,x_m^p)$. Since the
 Frobenius map is injective, and thus, since $k$ is finite, surjective, it follows that
 $\sum_{r=1}^t \lambda_r v_r\in T(k)$ and for each $r$, $\lend(v_r)=0$ and for each $x\in X$ that appears in
 $v_r$, $\deg_{x}(v_r)<q$. By Lemma \ref{lemma: field identities}, $f=0$.
\end{proof}

\begin{theorem}
If $p=2$, then $T(\unitgrass)=\set x_1^2-x_1^{2q},\com x_1,{x_2}\endset^T$, otherwise
 $T(\unitgrass)=\set x_1^{qp}-x_1^p,\com x_1,{x_2,x_3}\endset^T$.
\end{theorem} 

\begin{proof}

 If $p=2$, let $U=\set x_1^{2q}-x_1^2,\com x_1,{x_2}\endset^T$, while if
 $p>2$, let  $U=\set x_1^{qp}-x_1^p,\com x_1,{x_2,x_3}\endset^T$. By
 Lemma \ref{lemma: regev's identities}, $x_1^{qp}-x_1^p\in
 T(\unitgrass)$, and certainly $\com x_1,{x_2,x_3}\in T(\unitgrass)$, 
 while if $p=2$, then $\unitgrass$ is commutative and so $\com
 x_1,{x_2}\in T(\unitgrass)$ in this case. Thus $U\subseteq
 T(\unitgrass)$. Note that $T^{(3)}\subseteq \set \com
 x_1,{x_2}\endset^T$, so $T^{(3)}\subseteq U$ in every case. Suppose
 that $U\ne T(\unitgrass)$, and let $f\in T(\unitgrass)-U$. Since
 $x_1^{qp}-x_1^p\in U$, we may assume that for any variable $x$, if
 $x^\gamma$ is a factor of a summand of $f$, then $\gamma<qp$. As well,
 since $T^{(3)}\subseteq U$, we may assume that  $f=\sum_{i=1}^k
 \lambda_i v_i$ for $\lambda_i\in k^*$ and $v_i\in\siderovset$.
 Furthermore,  by Lemma \ref{lemma: p powers commute} and the fact that
 $T^{(3)}\subseteq U$, if $x\in X$ and $x^\gamma$ is a factor of $v_i$,
 with $\gamma=p\delta+\epsilon$ where $0\le\epsilon<p$, then we may move
 $x^{p\delta}$ to the front of $v_i$. Thus we may assume that
 $f=\sum_{i=1}^s f_i u_i$, where for each $i$, $f_i\ne 0$ is of the form
 described in Corollary \ref{corollary: needed for main identity
 theorem}, and either $u_1=1$ or else $u_i\in\boundedsiderovset$. By
 Corollary \ref{corollary: needed for main identity theorem}, not every
 $u_i$ is equal to 1. Let us represent the sum of the terms for which
 $u_i=1$ by $f_0$, and assume that the other terms have been labelled so
 that $f=f_0+ \sum_{i=1}^t f_i u_i$, where $u_1>u_2>\cdots >u_t$. By 
 Corollary \ref{corollary: needed for main identity theorem},  since
 $f_1\ne 0$, $f_1\notin T(\unitgrass)$. Suppose that the variables
 appearing in $f_1$ are $x_{i_1},x_{i_2},\ldots,x_{i_w}$. Then there
 exist $g_1,g_2,\ldots,g_w\in\unitgrass$ such that
 $f_1(g_1,g_2,\ldots,g_w)\ne 0$. Now, for any $g\in \unitgrass$, 
 $g^p=\lambda^p$, where $\lambda\in k$, $c\in C$, and $h\in H$ are such
 that $g=\lambda+c+h$. Since each variable of $f_1$ has degree a
 multiple of $p$, it follows that there are
 $\lambda_{i_1},\lambda_{i_2},\ldots,\lambda_{i_w}\in k$ such that
 $f_1(\lambda_{i_1},\lambda_{i_2},\ldots,\lambda_{i_w})\ne 0$. For any
 $i\ge1$, $i\notin \set i_1,i_2,\ldots, i_w\endset$, let $\lambda_i=0$.
 Now apply Proposition \ref{proposition: main evaluation fact} to obtain
 a homomorphism $\varphi\from \konex \to \unitgrass$ such that 
 $\varphi(f_1)=f_1(\lambda_{i_1},\ldots,\lambda_{i_w})\ne 0$ and 
 $\varphi(u_1)\ne 0$, while $\wt(\varphi(u_1))>\wt(\varphi(u_i))$ for
 all $i>1$. Since $\varphi(f_0)\in k$, it follows that  $\varphi(f)\ne
 0$, contradicting the fact that $f\in T(G)$. Thus $T(G)-U=\nullset$,
 and so $T(G)\subseteq U$,  as required.
\end{proof}


\section{The central polynomials of the infinite dimensional unitary Grassmann
           algebra over a finite field}  

 Recall that the $T$-space of
 $\konex$ that is generated by the set of all central polynomials of
 $\unitgrass$ is denoted by $CP(\unitgrass)$. By Lemma \ref{lemma: useful} (iv),
 $g^p\in C_{\unitgrass}$ for all $g\in \unitgrass$, and so $x^p\in CP(\unitgrass)$.
 
 When $p=2$, then $\unitgrass$ is commutative and so $CP(\unitgrass)=\konex$.
 Thus for the sequel, we assume that $p>2$.

\begin{proposition}\label{proposition: evaluation for central polynomials}
  Let $u\in\boundedsiderovset$ be such that $\beg(u)>0$, and let $t$ be such that $x_t$ appears in
  $\beg(u)$. Set $m=2\deg(u)-2\lend(u)-1$, and for each $i\ge1$, let $\lambda_i\in k$.
  Then there exists a homomorphism $\varphi\from \konex\to \finiteunitgrass{m}$ such that the following hold.
   \begin{list}{(\roman{parts})}{\usecounter{parts}}
   \item For each index $i$, $proj_k(\varphi(x_i))=\lambda_i$.
   \item $\displaystyle\dom(\varphi(u))=\lambda_t 2^{\lend(u)}\mkern-25mu\prod_{x\text{ in }\beg(u)}\mkern-25mu \deg_{x}(u)\mkern2mu!
    \mkern-25mu\prod_{x\text{ in }\endof(u)}\mkern-25mu (\deg_{x}(u)-1)\mkern2mu!\mkern5mu\prod_{i=1}^{m}e_i$.
   \item For any $v\in\boundedsiderovset$ with $u>v$ and either $x_t$ appears in $\beg(v)$ or else $x_t$ appears with degree $p$
  in $\endof(v)$, $m=\wt(\varphi(u))>\wt(\varphi(v))$.
  \end{list}
\end{proposition}

\begin{proof}

 The homomorphism $\varphi$ is determined by the following assignments.
 First, any variable $x_i\in X$ that does not appear in $u$ is mapped to
 $\lambda_i$. Then choose an index offset $N_t$, let
 $\alpha=\deg_{x_t}(u)$, set  $E_t=\set e_{N+\epsilon}\mid
 \epsilon=1,2,\ldots,2\alpha-1\endset\subseteq E$, and map $x_t$ to
 $\lambda_t+e_{N+2\alpha-1}+\sum_{\epsilon=1}^{\alpha-1}
 e_{N+2\epsilon-1}e_{N+2\epsilon}$. Next,  for any variable $x\ne x_t$
 that appears in $\beg(u)$, choose an index offset $N=N_x$, let
 $E_x=\set e_{N+\epsilon}\mid
 \epsilon=1,2,\ldots,2\deg_{x}(u)\endset\subseteq E$, and map  $x$ to
 $\lambda_x+\sum_{\epsilon=1}^{\deg_{x}(u)}
 e_{N+2\epsilon-1}e_{N+2\epsilon}$. Finally, for  any variable $x$ that
 appears in $\endof(u)$, choose an index offset $N=N_x$,  $E_x=\set
 e_{N+\epsilon}\mid \epsilon=1,2,\ldots,2\deg_{x}(u)-1\endset\subseteq
 E$, and map $x$ to
 $\lambda_x+e_{N+2\deg_{x}(u)-1}+\sum_{\epsilon=1}^{\deg_{x}(u)-1}
 e_{N+2\epsilon-1}e_{N+2\epsilon}$.  We observe that (i) is satisfied by
 this assignment. The offsets $N_x$ are chosen so that $x\ne y$ implies
 that  $E_x\cap E_y=\nullset$ and $\bigcup_{x\text{ appears in }u}
 E_x=\set e_i\mid i=1,2,\ldots,m\endset$.
 
 Recall that for $g_1,g_2\in\unitgrass$, $\dom(g_1g_2)=\dom(g_1)\dom(g_2)$
 if $s(g_1)\cap s(g_2)=\nullset$. In
 particular, since $u\in\boundedsiderovset$ (where the case of $u$ with
 empty end is just a simplification of the following argument), $u$ is
 of the form  
 $$
 \prod_{r=1}^tx_{i_r}^{\alpha_r}\prod_{r=1}^s [\mkern1mu
 x_{j_{2r-1}},{x_{j_{2r}}}]x_{j_{2r-1}}^{\beta_{2r-1}}x_{j_{2r}}^{\beta_{2r}},
 $$
 where for each $r=1,2,\ldots,t$, $1\le \alpha_{r}\le p-1$ and for each
 $r=1,2,\ldots,2s$, $0\le \beta_{r}\le p$, so $\dom(\varphi(u))$ will be
 the product of $\dom(\varphi(x_{i_r})^{\alpha_r}))$, $r=1,2,\ldots,t$,
 and  $$
  \dom([\mkern1mu \varphi(x_{j_{2r-1}}),\varphi({x_{j_{2r}}})]\varphi(x_{j_{2r-1}})^{\beta_{2r-1}}\varphi(x_{j_{2r}})^{\beta_{2r}})
 $$
 for $r=1,2,\ldots,s$. 

 We now apply Lemma \ref{technical result} to evaluate the dominant part of $\varphi(u)$, where for convenience,
 we shall let $N_{i}$ denote $N_x$ where $x=x_{i}$. We find that
 \begin{align*}
 \dom(\varphi(u))&=\prod_{r=1}^t 
  \biggl(\alpha_r!\prod_{\epsilon=1}^{2\alpha_{r}}e_{N_{i_r}+\epsilon}\biggr)\mkern4mu 2^{\lend(u)}\\
  &\hskip50pt\prod_{r=1}^s  \biggl( \beta_{2r-1}!
  \prod_{\epsilon=1}^{2\beta_{2r-1}+1} \mkern-20mu e_{N_{j_{2r-1}}+\epsilon}\mkern6mu  
  \beta_{2r}!\,\prod_{\epsilon=1}^{2\beta_{2r}+1}e_{N_{j_{2r}}+\epsilon}\biggr) \\
  &=2^{\lend(u)}\prod_{r=1}^t\alpha_r!\prod_{r=1}^{2s}\beta_{r}!\mkern20mu \prod_{\epsilon=1}^m e_\epsilon
 \end{align*} 
 and so (ii) holds. 
 
 Finally, suppose that $v\in \boundedsiderovset$, with $u>v$ and either
 $x_t$ appears in $\beg(v)$ or else $x_t$ appears in $\endof(v)$ with
 degree $p$. If $\varphi(v)=0$, then $\wt(\varphi(v))=0$, while
 $\wt(\varphi(u))=m>0$. Thus we may suppose that $\varphi(v)\ne 0$. If a
 variable other than $x_t$ appears in $\beg(u)$ and in $\endof(v)$, then
 $\varphi(v)=0$, so we may assume that other than $x_t$, every variable
 that appears in $\beg(u)$ and in $v$ then appears in $\beg(v)$. Since
 conditions (ii) and (iv) of the definition of the Siderov ordering
 imply that there is a variable that appears in $\beg(u)$ and in
 $\endof(v)$, we see that $u>v$ must hold by virtue of conditions (i) or
 (iii) of the definition, or else by conditions (ii) or (iv).  We shall
 prove that it is not possible for $u>v$ to hold by virtue of conditions
 (ii) or (iv), but first, let us consider the situation when $u>v$ by
 virtue of conditions (i) or (iii). We now calculate the weight of the
 dominant part of $\varphi(v)$. The contribution to  the weight of the
 dominant part of $\varphi(v)$ by $x\in X$ that appears in $v$ but not
 in $u$ is 0. By Lemma \ref{technical result}, the contribution to the
 weight of the dominant part of $\varphi(v)$ by $x\in X$ that appears in
 both $u$ and $v$  is $\min\set 2\deg_{x}(u),2\deg_x(v)\endset$ if $x$
 appears in $\beg(u)$ (so the contribution of $x$ to the dominant part
 of $\varphi(u)$ has weight $2\deg_x(u)$), while it is $\min\set
 2\deg_x(u)-1, 2\deg_x(v)-1\endset$ if  $x$ appears in $\endof(u)$ and
 $\endof(v)$ (so the contribution of $x$ to the dominant part of
 $\varphi(u)$ has weight $2\deg_x(u)-1$), and it is $\min\set
 2\deg_x(u)-1, 2\deg_x(v)\endset$ if $x$ appears in $\endof(u)$ and in
 $\beg(v)$ (so the contribution of $x$ to the dominant part of
 $\varphi(u)$ has weight $2\deg_x(u)$). Since both conditions (i) and
 (iii) imply that there is a variable $x$ such that
 $\deg_{x}(u)>\deg_{x}(v)\ge0$, it follows that
 $\wt(\varphi(u))>\wt(\varphi(v))$. Now suppose that $u>v$ by virtue of
 condition (ii) or (iv), either of which implies that there is $x\in X$
 such that $x$ appears in $\beg(u)$ and in $\endof(v)$. If $x\ne x_t$,
 then  $\varphi(v)=0$, so we may suppose that $x=x_t$. Recall that
 $\deg_{x_t}(v)=p$, while $\deg_{x_t}(u)\le p-1$. Since condition (iv)
 only applies when all variables have the same degree in both $u$ and
 $v$, we see that we must have $u>v$ by virtue of condition (ii). Thus
 $\deg(u)=\deg(v)$, but $\lend(v)>\lend(u)$. Now, every variable (other
 than $x_t$) that appears in $\beg(u)$ and also appears in $v$ appears
 in $\beg(v)$, and every variable that appears in $v$ also appears in
 $u$ (otherwise $\varphi(v)=0$), so the only variables that can appear
 in $\endof(v)$ are those in $\endof(u)$, contradicting the fact that
 $\lend(v)>\lend(u)$. This completes the proof that $u>v$ can't occur by
 virtue of conditions (ii) or (iv).
\end{proof} 
 
\begin{definition}\label{definition: def of S1}
  $$
   S_1=\set \com x_1,{x_2},x_1^p\endset^S+\set x_{1}^p
    \prod_{i=1}^t\com x_{2i},{x_{2i+1}}x_{2i}^{p-1}x_{2i+1}^{p-1}\mid t\ge1\endset^S.
  $$    
\end{definition}

The content of following lemma is known (for example, see Grishin and Tsybulya, \cite{GT}).

\begin{lemma}\label{lemma: for unitary}
 Let $u\in \siderovset$ with $\lbeg{u}>0$. 
 If $\deg_{x}(u)\cong0\mod{p}$ for every $x\in X$ that appears in the beginning
 of $u$, then $u\in S_1+T^{(3)}$.
\end{lemma}

\begin{theorem}\label{theorem: central in unitary}
 $CP(\unitgrass)=S_1+T(\unitgrass)$.
\end{theorem}

\begin{proof}
 Let $U_1=S_1+T(\unitgrass)$. By Lemma \ref{lemma: p powers commute} and Lemma 1.1 (vi) of \cite{Ra},
 $U_1\subseteq CP(\unitgrass)$. Suppose that $CP(\unitgrass)-U_1\ne\nullset$, and let $f\in
 CP(\unitgrass)-U_1$. Recall that $f$ is said to be {\it  essential in its variables} if every 
 variable that appears in any monomial of $f$ appears in every monomial of $f$. It is well-known
 that every $T$-space is generated (as a $T$-space) by
 its essential elements, so we may assume that $f$ is essential, and that there is no essential
 element of $CP(\unitgrass)-U_1$ with fewer summands. We may further assume that the variables
 that appear in $f$ are $x_1,x_2,\ldots,x_n$. There exist $v_1,v_2,\ldots,v_l\in\siderovset$ and 
 scalars $\alpha_1,\alpha_2,\ldots, \alpha_l$ such that, modulo $T^{(3)}$, $f=\sum_{i=1}^l \alpha_iv_i$. 
 If any $v_i$ had $\deg_x(v_i)\cong 0\mod{p}$ for each $x\in X$ that appears in 
 $\beg(v_i)$, then $v_i\in U_1$ by Lemma \ref{lemma: for unitary}. Since this would contradict our 
 choice of $f$, it follows that for each $i$, there exists $x\in X$
 that appears in $\beg(v_i)$ with degree not a multiple of $p$. Thus, for each $i$, we may (modulo $T^{(3)}$)
 write $v_i$ in the form $f_iu_i$, where $f_i$ is a product of the form $\prod_{r=1}^t x_{i_r}^{p\delta_r}$, 
 with $i_1<i_2<\cdots i_t$, and $u_i\in\boundedsiderovset$ such that 
 $\set x_1,x_2,\ldots,x_n\endset=\set x_{i_1},\ldots,x_{i_t}\endset
 \cup\set x\in X\mid x\text{ appears in $u_i$}\endset$. Finally, since $x_i^p-x_i^{qp}\in T(\unitgrass)$,
 we may assume that each $f_i$ is a $p$-polynomial. We remark that $\set x_{i_1},\ldots,x_{i_t}\endset$
 and $\set x\in X\mid x\text{ appears in $u_i$}\endset$ need not be disjoint. Now, although $v_i\ne v_j$ if $i\ne j$,
 the same need not be true for the $u_i$'s. Therefore, we shall write $f=\sum_{i=1}^s f_iu_i$, where now
 $u_1>u_2>\cdots>u_s$ and for each $i$, $f_i$ is an essential $p$-polynomial. Choose an index $t$ such that
 $x_t$ appears in $\beg(u_1)$. Consider any $i>1$ such that $x_t$ appears in $\endof(u_i)$ with
 $\deg_{x_t}(u_i)<p$. By Corollary 2.2 of \cite{Ra} 
 (where the $T$-space $S$ referred to there is
 a subspace of $S_1$), there exists $u_i'\in \boundedsiderovset$ and $\lambda_i\in k$ such that $\lambda_i\ne 0$,
 $u_i\cong \lambda_iu_i'\mod{U_1}$, and $x_t$ appears in $\beg(u_i')$. In the expression for $f$, replace
 each such $u_i$ by the corresponding $\lambda_iu_i'$. By the choice of $f$, no cancellation can occur (although it is possible
 that $u_i'=u_j'$ for $i\ne j$). Relabelling as necessary, we now have an essential element 
 $\sum_{i=1}^q f_iu_i\in CP(\unitgrass)-U_1$ such that $u_1>u_2>\cdots > u_q$, $x_t$ appears in $\beg(u_1)$, 
 and for each $i$, $f_i$ is an essential $p$-polynomial and either $x_t$ appears in $\beg(u_i)$
 or else $x_t$ appears in $\endof(u_i)$ with degree $p$. Of all such polynomials, let $g$ denote one
 for which $u_1$ is least. Write $g=\sum_{x_t\text{ appears in $\beg(u_i)$}}f_iu_i + \sum_{x_t\text{ appears in $\endof(u_i)$}}f_iu_i$
 and observe that by choice of $g$, $\sum_{x_t\text{ appears in $\beg(u_i)$}}f_iu_i\ne 0$.
 If $f_1\in k$, then $f_1\ne0$. Otherwise, $f_1\notin k$, and then by Corollary \ref{corollary: needed for main identity theorem}, 
 $f_1\notin T(\unitgrass)$.In this case, suppose that the variables that appear in $f_1$ are $x_{i_1},x_{i_2},\ldots,x_{i_d}$.
 Since for any $g\in \unitgrass$, $g^p=proj_k(g)^p$, it follows that there are $\lambda_{i_1},\lambda_{i_2},\ldots,\lambda_{i_d}\in k$
 such that $f_1(\lambda_{i_1},\lambda_{i_2},\ldots,\lambda_{i_d})\ne 0$. Note that since $f_1$ is essential, $\lambda_{i_j}\ne 0$ for
 every $j$. For any $i\notin\set i_1,i_2,\ldots,i_d\endset$,  set $\lambda_i=1$. Note that $\lambda_t\ne 0$.  Now 
 apply Proposition \ref{proposition: evaluation for central polynomials} to $u_1$, $x_t$, and our selected $\lambda_i$'s to obtain 
 a homomorphism $\varphi\from \konex\to \finiteunitgrass{m}$ such that for $m=2\deg(u_1)-2\lend(u_1)-1$ the following hold.
   \begin{list}{(\roman{parts})}{\usecounter{parts}}
   \item For each index $i$, $proj_k(\varphi(x_i))=\lambda_i$.
  \item $\displaystyle\dom(\varphi(u_1))=\lambda_t 2^{\lend(u_1)}\mkern-25mu\prod_{x\text{ in }\beg(u_1)}\mkern-25mu \deg_{x}(u_1)\mkern2mu!
    \mkern-25mu\prod_{x\text{ in }\endof(u_1)}\mkern-25mu (\deg_{x}(u_1)-1)\mkern2mu!\mkern5mu\prod_{i=1}^{m}e_i$.
  \item For any $v\in\boundedsiderovset$ with $u_1>v$ and either $x_t$ appears in $\beg(v)$ or else $x_t$ appears with degree $p$
  in $\endof(v)$, $m=\wt(\varphi(u_1))>\wt(\varphi(v))$.
  \end{list}
 It follows that 
 $$
 \varphi(f)=\lambda\mkern5mu\prod_{i=1}^{m}e_i+
     \text{terms of lower weight}
 $$     
 where 
 $$
  \lambda=f_1(\lambda_{i_1},\ldots,\lambda_{i_d})\lambda_t 2^{\lend(u_1)}\mkern-25mu\prod_{x\text{ in }\beg(u_1)}\mkern-25mu \deg_{x}(u_1)\mkern2mu!
    \mkern-25mu\prod_{x\text{ in }\endof(u_1)}\mkern-25mu (\deg_{x}(u_1)-1)\mkern2mu!\ne 0.
 $$     
 Since $m$ is odd, we obtain that $\varphi(f)$ has nonzero odd part. But then $\varphi(f)\notin C_{\unitgrass}$, which contradicts the fact
 that $f\in CP(\unitgrass)$. Thus $CP(\unitgrass)-U_1=\nullset$.
\end{proof}

\end{document}